\documentclass[12pt]{amsart}
\usepackage{amsmath,amsthm,latexsym,amscd,amsbsy,amssymb,url}
\usepackage[all]{xypic}
 \relax

\chardef\bslash=`\\ 

\makeatletter
\def\verbatim{\interlinepenalty\@M \@verbatim
  \leftskip\@totalleftmargin\advance\leftskip2pc
  \frenchspacing\@vobeyspaces \@xverbatim}
\makeatother
\makeatletter
  \def\dgt@k{\dg@DX=-3 \dg@DY=2 \dg@SIZE=3}
\makeatother

\makeatletter
  \def\dgt@kk{\dg@DX=3 \dg@DY=-1 \dg@SIZE=3}%
\makeatother

\theoremstyle{plain}
\newtheorem{thm}{Theorem}[section]

\newtheorem{lem}[thm]{Lemma}
\newtheorem{pro}[thm]{Proposition}

\theoremstyle{definition}

\numberwithin{equation}{section}

\begin{document}


\title
{Valdivia compact groups are products}
\author{A. Chigogidze}
\address{Department of Mathematics and Statistics,
University of North Carolina at Greensboro,
383 Bryan Bldg, Greensboro, NC, 27402, USA}
\email{chigogidze@uncg.edu}
\keywords{Compact group, Valdivia compact, inverse spectrum}
\subjclass{Primary: 54D30; Secondary: 54C15.}


\begin{abstract}{It is shown that every Valdivia compact group is homeomorphic to a product of metrizable compacta.}
\end{abstract}

\maketitle
\markboth{A.~Chigogidze}
{Remarks on Valdivia compact spaces}

\section{Introduction}\label{S:intro}
A compact space is called Valdivia compact if it can be embedded into ${\mathbb R}^{A}$ for some $A$ in such a way that the image of this embedding is the closure of a subset of a $\Sigma$-product. These spaces have been studied by several authors (see \cite{kalenda}, \cite{kalenda1} for a nice survey of Valdivia compacta and related topics). A particular area of interest is understanding of the structure of Valdivia compact groups. There are compact connected Abelian groups which are not Valdivia compact \cite{ku}. On the other hand it was shown recently \cite{kubisgroups} that every Valdivia compact Abelian group is homeomorphic to a product of metrizable compacta. We extend this result to the noncommutative case (Theorem  \ref{T:main}). This  answers negatively question from \cite{kubisgroups} whether there exist Valdivia compact groups which are not products.

Our original proof was based on a recently announced result on preservation of the class of Valdivia compacta by retractions. I would like to thank the referee for pointing out that the proof of this statement contained an error and that the validity of the result itself is still open. The proof of Theorem \ref{T:main} has been revised and as presented below does not depend on the above mentioned statement.

\section{Preliminaries}\label{S:pre}
For the reader's convenience in this section we present some of the needed results.
\subsection{Directed sets}\label{SS:dirset}
Let $\kappa \geq \omega$. A
subset  $B$  of a partially ordered directed set $A$ is said to be ${\kappa}$-{\em
closed} in  $A$  if for each chain
$C {\subseteq} B$, with $|C| \leq \kappa$, we have  $\sup C \in B$
whenever the element
$\sup C$ exists in $A$ ($\sup C$ denotes the lowest upper bound of elements of $C$). 
A set  $A$  is said to be
${\kappa}$-{\em complete}
if for each chain
$B$  of
elements of $A$, with $|B| \leq \kappa$,
there exists an element  ${\sup}B$ in $A$.
A standard example of
a ${\kappa}$-complete set is the set  ${\exp}_{\kappa}A$
of all subsets of cardinality $\leq \kappa$ of any 
set $A$. 

The following statement presents an important property of $\kappa$-complete sets (\cite[Proposition 1.1.27]{chibook}). 

\begin{pro}\label{P:3.1.1}
Let  $\{ A_i \colon i \in I \}$, $|I| \leq \kappa$, be a 
collection of
$\kappa$-closed and cofinal subsets of a
$\kappa$-complete
set $A$. Then the
intersection $\cap \{ A_i \colon i \in I \}$
is also cofinal
and $\kappa$-closed in $A$.
\end{pro}

\subsection{Inverse Spectra}\label{SS:spectra}
All limit projections of inverse spectra considered below are surjective and all spaces are compact. Let $\kappa \geq \omega$. An inverse spectrum ${\mathcal S}_{X} = \{ X_{\alpha}, p_{\alpha}^{\beta}, A\}$ consisting of compact spaces is a $\kappa$-spectrum if:
\begin{itemize}
\item[(i)]
$w(X_{\alpha}) \leq \kappa$, $\alpha \in A$;
\item[(ii)]
The indexing set $A$ is $\kappa$-complete;
\item[(iii)]
${\mathcal S}_{X}$ is $\kappa$-continuous, i.e. for each chain $\{\alpha_{i}\colon i \in I\} \subseteq A$ with $|I|\leq \kappa$ and $\alpha = \sup\{ \alpha_{i} \colon i \in I\}$, the diagonal product $\triangle\{ p_{\alpha_{i}}^{\alpha} \colon i \in I\} \colon X_{\alpha} \to \lim\{ X_{\alpha_{i}}, p_{\alpha_{i}}^{\alpha_{j}}, I\}$ is a homeomorphism.
\end{itemize}

One of the main results concerning $\kappa$-spectra is the following result of \v{S}\v{c}epin (known as the \v{S}\v{c}epin's Spectral Theorem, see \cite[Theorem 1.3.4]{chibook}). 

\begin{thm}\label{T:spectraltheorem}
Let ${\mathcal S}_{X} = \{ X_{\alpha}, p_{\alpha}^{\beta}, A\}$ and ${\mathcal S}_{Y} = \{ Y_{\alpha}, q_{\alpha}^{\beta},A\}$ be two $\kappa$-spectra. Then for every map $f \colon \lim{\mathcal S}_{X} \to\lim{\mathcal S}_{Y}$ there exist a cofinal and $\kappa$-closed subset $B \subseteq A$ and maps $f_{\alpha} \colon X_{\alpha} \to Y_{\alpha}$, $\alpha \in B$, such that $f = \lim\{ f_{\alpha} \colon \alpha \in B\}$. If $f$ is a homeomorphism, then we may assume that each $f_{\alpha}$, $\alpha \in B$, is also a homeomorphism.
\end{thm}

\subsection{$\Sigma$-products}\label{SS:sigma}
If $C \subseteq B \subseteq A$, then $\pi_{B} \colon {\mathbb R}^{A} \to {\mathbb R}^{B}$ and $\pi_{C}^{B} \colon {\mathbb R}^{B} \to {\mathbb R}^{C}$ denote the corresponding projections. Similarly by $i_{B} \colon {\mathbb R}^{B} \to {\mathbb R}^{A}$ and $i_{C}^{B} \colon {\mathbb R}^{C} \to {\mathbb R}^{B}$ we denote sections of $\pi_{B}$ and $\pi_{C}^{B}$ defined as follows:
\[ i_{B}(\{ x_{t} \colon t \in B\}) = \left(\{ x_{t} \colon t \in B\} ,\{ 0_{t} \colon t \in A\setminus B\}\right) \]

\noindent and

\[ i_{C}^{B}(\{ x_{t} \colon t \in C\}) = \left(\{ x_{t} \colon t \in C\} ,\{ 0_{t} \colon t \in B\setminus C\}\right) .\]

Below we consider the $\Sigma$-product of real lines, which is the subspace

\[ \Sigma (A) = \left\{ \{ x_{t} \colon  t \in A\} \in {\mathbb R}^{A} \colon |\{ t \in A \colon x_{t}\neq 0 \}|\leq \omega \}\right\} \]

\noindent of the product ${\mathbb R}^{A}$.

For each $B \subseteq A$, $\Sigma(B)$ is identified with the subspace $\Sigma(A) \cap i_{B}({\mathbb R}^{B})$ of $\Sigma(A)$. Note that $\Sigma(A) = \bigcup\{ i_{B}({\mathbb R}^{B}) \colon B \in \exp_{\omega}A\}$. 

Retractions $r_{B} = i_{B}\pi_{B} \colon \Sigma(A) \to \Sigma(B)$, $B \subseteq A$, play  an important role below. We refer the reader to \cite{gulko2} for relevant definitions and results. For a given closed subset $F \subseteq \Sigma(A)$ we are interested in finding subsets $B \subseteq A$ such that $r_{B}(F) \subseteq F$ or, equivalently, $F \cap \Sigma(B) = r_{B}(F)$. Note that if $Y = \operatorname{cl}_{{\mathbb R}^{A}}F$ is compact, then $r_{B}(F) \subseteq F$ if and only if $r_{B}(Y) \subseteq Y$ or, equivalently, $\operatorname{cl}(Y \cap \Sigma(B)) = r_{B}(Y)$. In this case $B$ is called an $F$-good (or $Y$-good) subset of $A$.  It turns out (see \cite[Lemma 1]{gulko1}, \cite[Lemma 1.2]{argyros}) that there are many $Y$-good subsets . Precise statement we use below is recorded in the following lemma.

\begin{lem}\label{L:good}
Let $\omega \leq \kappa < \tau$, $|A| = \tau$ and $Y = \operatorname{cl}_{{\mathbb R}^{A}}(Y \cap \Sigma(A))$ be compact. Then 
\begin{itemize}
\item[(a)]
the set ${\mathcal A}_{\kappa}^{V}$ of $Y$-good subsets of $A$ of cardinality $\leq \kappa$ is cofinal and $\kappa$-closed in $\exp_{\kappa}A$;
\item[(b)]
union of an increasing collection of $Y$-good subsets of $A$ is again a $Y$-good subset of $A$;
\item[(c)]
if $B$ is a $Y$-good subset of $A$, then $r_{B}|Y \colon Y \to Y_{B}$ is a retraction.
\end{itemize} 
\end{lem}



\section{Proof of the Main Result}\label{S:main}
In this section we prove our main result - Theorem \ref{T:main}. First we need the following statement.

\begin{lem}\label{L:group}
Let $\omega \leq \kappa < \tau$, $|A| = \tau$ and a compact group $X$ be embedded into  ${\mathbb R}^{A}$. Then the set $\exp_{\kappa}A$ contains a cofinal and $\kappa$-closed subset ${\mathcal A}_{\kappa}^{G}$ such that $\pi_{B}(X)$ is a topological group and the projection $\pi_{B}|X \colon X \to \pi_{B}(X)$ is a homomorphism for each $B\in {\mathcal A}_{\kappa}^{G}$.
\end{lem}
\begin{proof}
Let $\lambda \colon X \times X \to X$ and $\mu \colon X\to X$ denote group operations in $X$, i.e. $\lambda (x,y) = x\cdot y$ and $\mu(x) = x^{-1}$, $x,y \in X$. Let also $X_{B} = \pi_{B}(X)$, $p_{B} = \pi_{B}|X$ and $p_{C}^{B} = \pi_{C}^{B}|X_{B}$ for any $C, B \subseteq A$ with $C \subseteq B$. 

Consider the standard $\kappa$-spectrum ${\mathcal S}_{X} = \{ X_{B}, p_{C}^{B}, \exp_{\kappa}A\}$ and note that $X = \lim {\mathcal S}_{X}$. By Theorem \ref{T:spectraltheorem}, applied to the map $\mu \colon X \to X$ and the $\kappa$-spectrum ${\mathcal S}_{X}$, there exist a cofinal and $\kappa$-closed subset ${\mathcal A}_{\kappa}^{\mu}$ of $\exp_{\kappa}A$ and maps $\mu_{B} \colon X_{B} \to X_{B}$, $B \in {\mathcal A}_{\kappa}^{\mu}$, such that $p_{B}\mu = \mu_{B}p_{B}$ for each $B \in {\mathcal A}_{\kappa}^{\mu}$ and $\mu = \lim \{ \mu_{B} \colon B \in {\mathcal A}_{\kappa}^{\mu}\}$.

Next consider the $\kappa$-spectrum ${\mathcal S}_{X} \times {\mathcal S}_{X}  = \{ X_{B}\times X_{B}, p_{C}^{B}\times p_{C}^{B}, \exp_{\kappa}A\}$. Obviously, $X \times X = \lim {\mathcal S}_{X}\times {\mathcal S}_{X}$. Applying Theorem \ref{T:spectraltheorem} to the map $\lambda \colon X \times X \to X$ and to the spectra ${\mathcal S}_{X}\times {\mathcal S}_{X}$ and ${\mathcal S}_{X}$, we conclude that there exist a cofinal and $\kappa$-closed subset ${\mathcal A}_{\kappa}^{\lambda}$ of $\exp_{\kappa}A$ and maps $\lambda_{B} \colon X_{B} \times X_{B} \to X_{B}$, $B \in {\mathcal A}_{\kappa}^{\lambda}$ such that $p_{B}\lambda = \lambda_{B}(p_{B}\times p_{B})$ for each $B \in {\mathcal A}_{\kappa}^{\lambda}$ and $\lambda = \lim \{ \lambda_{B} \colon B \in {\mathcal A}_{\kappa}^{\lambda}\}$

By Proposition \ref{P:3.1.1}, the intersection ${\mathcal A}^{G}_{\kappa} = {\mathcal A}_{\kappa}^{\mu} \cap {\mathcal A}_{\kappa}^{\lambda}$ is still cofinal and $\kappa$-closed in $\exp_{\kappa}A$. Note that for each $B \in {\mathcal A}_{\kappa}^{G}$ we have two maps $\lambda_{B} \colon X_{B} \times X_{B} \to X_{B}$ and $\mu_{B} \colon X_{B} \to X_{B}$. These maps define a group structure on $X_{B}$ as follows: $x_{B}\cdot y_{B} = \lambda_{B}(x_{B},y_{B})$ and $x_{B}^{-1} = \mu_{B}(x_{B})$, $x_{B}, y_{B} \in X_{B}$. The unit element in $X_{B}$ is defined as follows: $e_{B} = p_{B}(e)$, where $e$ is the unit in $X$. It is easy to see that the projection $p_{B}\colon X \to X_{B}$ becomes a group homomorphism. Indeed, for $x, y \in X$ we have
\[ p_{B}(x\cdot y) = p_{B}\lambda(x,y) = \lambda_{B}(p_{B}\times p_{B})(x,y) = \lambda_{B}(p_{B}(x),p_{B}(y)) = p_{B}(x)\cdot p_{B}(y).\]\end{proof}

\begin{lem}\label{L:retract}
Let $X$ be a topological group which is a retract of a Valdivia compact. Then there exists a well-ordered continuous spectrum ${\mathcal S}_{X} = \{ X_{\alpha}, p_{\alpha}^{\alpha +1},\tau\}$ such that
\begin{enumerate}
\item 
$\tau = w(X)$;
\item 
$w(X_{\alpha}) \leq |\omega| \cdot |\alpha|$, $\alpha < \tau$;
\item
$X_{\alpha}$ is a topological group and a retract of a Valdivia compact, $\alpha < \tau$;
\item
The limit projection $p_{\alpha} \colon X \to X_{\alpha}$ is a topological homomorphism and a retraction, $\alpha < \tau$.
\end{enumerate}
\end{lem}
\begin{proof}
Let $Y$ be a Valdivia compact of suitably embedded into ${\mathbb R}^{A}$ with $|A| = \tau > \omega$. Suppose also that $s \colon Y \to X$ is a retraction and $X$ is a topological group. Without loss of generality we may assume that $w(X) = \tau$. 

Let $\omega \leq \kappa < \tau$. By Lemma \ref{L:good}(a),  the collection ${\mathcal S}_{Y}^{\kappa} = \{ Y_{B}, q_{C}^{B}, C,B \in {\mathcal A}_{\kappa}^{V}\}$, where $q_{C}^{B} = r_{C}^{B}|Y_{B}$, whenever $C, B \in {\mathcal A}_{\kappa}^{V}$ and $C \subseteq B$, forms a $\kappa$-spectrum. Also consider the $\kappa$-spectrum ${\mathcal S}_{X}^{\kappa} = \{ X_{B}, p_{C}^{B}, C,B \in {\mathcal A}_{\kappa}^{V}\}$, where $X_{B} = r_{B}(X)$ and $p_{C}^{B} = r_{C}^{B}|X_{B}$, whenever $C, B \in {\mathcal A}_{\kappa}^{V}$ and $C \subseteq B$. Clearly, $Y = \lim{\mathcal S}_{Y}^{\kappa}$ and $X = \lim{\mathcal S}_{X}^{\kappa}$. Note that projections $q_{C}^{B} \colon Y_{B} \to Y_{C}$ and $q_{B} \colon Y \to Y_{B}$ of the spectrum ${\mathcal S}_{Y}^{\kappa}$ are retractions (Lemma \ref{L:good}(c)). Theorem \ref{T:spectraltheorem}, applied to the $\kappa$-spectra ${\mathcal S}_{Y}^{\kappa}$ and ${\mathcal S}_{X}^{\kappa}$ and to the retraction $s \colon Y \to X$, guarantees existence of a $\kappa$-closed and cofinal subset ${\mathcal A}_{s}^{\kappa}$ of ${\mathcal A}_{\kappa}^{V}$ such that for each $B \in {\mathcal A}_{s}$ there is a map $s_{B} \colon Y_{B} \to X_{B}$ satisfying the equality $p_{B}s = s_{B}q_{B}$. 

Next observe that $s_{B} \colon Y_{B} \to X_{B}$, $B \in {\mathcal A}_{s}^{\kappa}$, is a retraction. Indeed, let $x_{B} \in X_{B}$ and $x \in X$ is such that $p_{B}(x) = x_{B}$. Then we have

\[  s_{B}(x_{B}) = s_{B}(p_{B}(x)) = s_{B}(q_{B}(x) = p_{B}(s(x)) = p_{B}(x) = x_{B}.\]

Next note that since in the diagram ($B \in {\mathcal A}_{s}^{\kappa}$)

 \[
        \xymatrix{
            X  \ar_{p_{B}}[d] & Y \ar_{s}[l] \ar^{q_{B}}[d]\\
            X_{B}  & Y_{B} \ar_{s_{B}}[l] \\
        }
      \]

\noindent maps $s$, $s_{B}$ and $q_{B}$ are retractions, it follows that $p_{B}$ is also a retraction.

By Lemma \ref{L:group}, we may assume without loss of generality that there exists a $\kappa$-closed and cofinal subset ${\mathcal A}^{\kappa}$ of ${\mathcal A}_{s}^{\kappa} \cap {\mathcal A}_{\kappa}^{G}$ such that for each $B \in {\mathcal A}^{\kappa}$ compactum $X_{B}$ is a topological group and the projection $p_{B} \colon X \to X_{B}$ is a topological homomorphism.

Thus, for any $\kappa$ with $\omega \leq \kappa < \tau$ we have a $\kappa$-closed and cofinal subset ${\mathcal A}^{\kappa}$ of $\exp_{\kappa}A$ such that the following conditions are satisfied for each $B \in {\mathcal A}^{\kappa}$:
\begin{itemize}
\item[(i)]
$Y_{B}$ is a Valdivia compact
\item[(ii)] 
$q_{B} \colon Y \to Y_{B}$ is a retraction
\item[(iii)]
There exists a retraction $s_{B} \colon Y_{B} \to X_{B}$ such that $p_{B}s = s_{B}q_{B}$. 
\item[(iv)] 
$p_{B} \colon X \to X_{B}$ is a retraction 
\item[(v)] 
$X_{B}$ is a topological group
\item[(vi)] 
$p_{B} \colon  X \to X_{B}$ is a topological homomorphism.
\end{itemize}

Let us now define subsets $A_{\alpha}$, $\alpha < \tau$, of $A$. Fix a well-ordering $\{ B_{\alpha} \colon \alpha < \tau\}$ of the set ${\mathcal A}^{\kappa}$, let $A_{0} = B_{0}$, $s_{0} = s_{A_{0}}$ and assume that the sets $A_{\beta}$ and retractions $s_{\beta} \colon Y_{A_{\beta}} \to X_{A_{\beta}}$ have been constructed for each $\beta < \alpha$, where $\alpha < \tau$, in such a way that the following conditions are satisfied: 

\begin{itemize}
\item[(a)$_{\beta}$]
$|A_{\beta}| \leq |\omega| \cdot |\beta|$
\item[(b)$_{\beta}$]
$A_{\gamma}\subseteq A_{\beta}$ whenever $\gamma \leq \beta < \alpha$
\item[(c)$_{\beta}$] 
$A_{\beta} = \cup\{ A_{\gamma} \colon \gamma < \beta\}$ whenever $\beta$ is a limit ordinal.
\item[(d)$_{\beta}$]
$A_{\beta +1} \supseteq A_{\beta} \cup B_{\beta +1}$, whenever $\beta +1 < \alpha$.
\item[(e)$_{\beta}$]
There exists a retraction $s_{\beta} \colon Y_{A_{\beta}} \to X_{A_{\beta}}$ such that $p_{A_{\beta}}s = s_{\beta}q_{A_{\beta}}$
\item[(f)$_{\beta}$] 
$X_{A_{\beta}}$ is a topological group
\item[(g)$_{\beta}$] 
$p_{A_{\beta}} \colon X \to X_{A_{\beta}}$ is a retraction and a topological homomorphism
\end{itemize}

If $\alpha$ is a limit ordinal, then let $A_{\alpha} = \cup\{ A_{\beta} \colon \beta < \alpha\}$. Note that for each $\gamma, \beta$ with $\gamma < \beta < \alpha$ the following diagram

 \[
        \xymatrix{
            X_{A_{\beta}}  \ar_{p_{A_{\gamma}}^{A_{\beta}}}[d] & Y_{A_{\beta}} \ar_{s_{\beta}}[l] \ar^{q_{A_{\gamma}}^{A_{\beta}}}[d]\\
            X_{A_{\gamma}}  & Y_{A_{\gamma}} \ar_{s_{\gamma}}[l] \\
        }
      \]

\noindent is commutative. Note also that $X_{A_{\alpha}} = \lim\{ X_{A_{\beta}}, p_{A_{\gamma}}^{A_{\beta}}, \alpha\}$ and $Y_{A_{\alpha}} = \lim\{ Y_{A_{\beta}}, q_{A_{\gamma}}^{A_{\beta}}, \alpha\}$. Consequently, collection $\{ s_{\beta} \colon \beta < \alpha\}$ defines a map $s_{\alpha} \colon Y_{A_{\alpha}} \to X_{A_{\alpha}}$ such that $p_{A_{\beta}}s_{\alpha} = s_{\beta}q_{A_{\beta}}$ for each $\beta < \alpha$. Since $s_{\beta}$ is a retraction, for any $x \in X_{A_{\alpha}}$ we have

\[ p_{A_{\beta}}(s_{\alpha}(x)) = s_{\beta}(q_{A_{\beta}}(x)) = s_{\beta}(p_{A_{\beta}}(x)) = p_{A_{\beta}}(x) ,\]

\noindent which shows that $s_{\alpha}(x) = x$. This verifies condition (e)$_{\alpha}$. By (b)$_{\beta}$ and Lemma \ref{L:good}(b), $A_{\alpha}$ is a $Y$-good subset of $A$. Consequently, by Lemma \ref{L:good}(c), $q_{A_{\alpha}} \colon Y \to Y_{A_{\alpha}}$ is a retraction. But then it is easy to conclude that $p_{A_{\alpha}} \colon X \to X_{A_{\alpha}}$ is also a retraction.

By conditions (f)$_{\beta}$ and (g)$_{\beta}$, it follows that for $\beta < \alpha$, $X_{A_{\alpha}} = \lim\{ X_{A_{\beta}}, p_{A_{\gamma}}^{A_{\beta}}, \alpha\}$ is also a topological group with naturally defined group operations and the projection $p_{A_{\alpha}} \colon X \to X_{A_{\alpha}}$ is a topological homomorphism. 

Next consider the case $\alpha = \beta +1$. Since, by condition (a)$_{\beta}$, $|A_{\beta}| \leq |\omega| \cdot |\beta|$ and since ${\mathcal A}^{\kappa}$ (with $\kappa = |\omega| \cdot |\beta|$) is cofinal in $\exp_{\kappa}A$, there exists $A_{\alpha} \in {\mathcal A}^{\kappa}$ such that $A_{\alpha} \supseteq A_{\beta} \cup B_{\alpha}$. Then there exists a retraction $s_{\alpha} \colon Y_{A_{\alpha}} \to X_{A_{\alpha}}$ such that $p_{A_{\alpha}}s = s_{\alpha}q_{A_{\alpha}}$, $p_{A_{\alpha}} \colon X \to X_{A_{\alpha}}$ is a retraction, $X_{A_{\alpha}}$ is a topological group and $p_{A_{\alpha}}$ is a topological homomorphism. 

This completes the construction. Let $X_{\alpha} = X_{A_{\alpha}}$ and $p_{\alpha} = p_{A_{\alpha}}$ for each $\alpha < \tau$. It only remains to note that $X = \lim\{ X_{\alpha}, p_{\alpha}^{\alpha +1},\tau\}$ and that all the required conditions (1)--(4) are satisfied.\end{proof}

\begin{thm}\label{T:main}
Let $X$ be a compact group. Then the following conditions are equivalent:
\begin{enumerate}
\item
$X$ is a Valdivia compact;
\item
$X$ is a retract of a Valdivia compact;
\item
$X$ is homeomorphic to a product of metrizable compacta;
\end{enumerate}
\end{thm}
\begin{proof}
Implications (1)$\Longrightarrow$ (2) and (3) $\Longrightarrow$ (1) are trivial. Let us prove imlication (2)$\Longrightarrow$ (3).
We proceed by induction. For metrizable groups which are retracts of Valdivia compacta there is nothing to prove. Let $X$ denote a topological group which is a retract of a Validivia compact. Suppose that the statement is true in cases when $w(X) < \tau$, $\tau > \omega$, and consider $X$ of weight $w(X) = \tau$. 
Let ${\mathcal S}_{X} = \{ X_{\alpha}, p_{\alpha}^{\alpha +1},\omega_{0} \leq \alpha <\tau\}$ be a spectrum supplied by Lemma \ref{L:retract}. Note that by property (4) in Lemma \ref{L:retract} each limit projection $p_{\alpha} \colon X \to X_{\alpha}$ is a retraction and a group homomorphism. This implies that each short projection $p_{\alpha}^{\alpha +1} \colon X_{\alpha +1} \to X_{\alpha}$ is a group homomorphism (and obviously a retraction). Consequently, $X_{\alpha +1}$ is homeomorphic to the product $X_{\alpha} \times \operatorname{ker}(p_{\alpha}^{\alpha +1})$ (see \cite[Lemma 0.2]{chi}) (corresponding homeomorphism is provided by the formula $h(x) = (p_{\alpha}^{\alpha +1}(x), x\cdot \left( p_{\alpha}^{\alpha+1}(x)\right)^{-1})$, $x \in X_{\alpha +1}$). Thus $X$ is homeomorphic to the product $X_{\omega_{0}}\times \prod\{ \operatorname{ker}(p_{\alpha}^{\alpha+1}) \colon \omega_{0}\leq \alpha < \tau\}$. Finally note that $\operatorname{ker}(p_{\alpha}^{\alpha +1})$ is a retract of $X_{\alpha +1}$ and $X_{\alpha +1}$, by Lemma \ref{L:retract}(3), is a retract of a Valdivia compact. Therefore $\operatorname{ker}(p_{\alpha}^{\alpha +1})$ is a retract of a Valdivia compact. By Lemma \ref{L:retract}(2), $w(\operatorname{ker}(p_{\alpha}^{\alpha +1})) \leq w(X_{\alpha +1}) \leq  |\omega| \cdot |\alpha +1| < \tau$. Then the inductive assumption guarantees that each $\operatorname{ker}(p_{\alpha}^{\alpha +1})$ is homeomorphic to a product of metrizable compact groups. Proof is completed.
\end{proof}

\end{document}